\documentclass{amsart}
\usepackage{palatino,hhline, bbm}
\usepackage{arydshln} 
\usepackage{amsthm, amsmath}
\usepackage{amssymb} 
\usepackage{amscd}
\usepackage{mathrsfs}
\renewcommand{\mathcal}{\mathscr}
\usepackage[hmargin=3cm, vmargin=2.5cm]{geometry}
\usepackage[breaklinks,bookmarksopen,bookmarksnumbered]{hyperref}
\usepackage[all]{xy}

\theoremstyle{plain}
\newtheorem*{thm*}{Theorem}

\theoremstyle{remark}

\newcommand\pr{\noindent\textit{Proof} : }

\newcommand\rond{\kern 1pt{\scriptstyle\circ}\kern 1pt}

\newcommand\End{\operatorname{End}}

\newcommand\Hom{\operatorname{Hom}}

\newcommand\pp{_{}^{\scriptscriptstyle\bullet}}

\renewcommand\P{\mathbb{P}}

\renewcommand\O{\mathcal{O}}

\newcommand\iso{\vbox{\hbox to .8cm{\hfill{$\scriptstyle\sim$}\hfill}
\nointerlineskip\hbox to .8cm{{\hfill$\longrightarrow $\hfill}} }}
\newcommand\bir{\vbox{\hbox to .8cm{\hfill{$\scriptstyle\sim$}\hfill}
\nointerlineskip\hbox to .8cm{{\hfill$\dasharrow $\hfill}} }}

\begin{document}
\title{Ulrich bundles on  surfaces with $q=p_g=0$}
\author[Arnaud Beauville]{Arnaud Beauville}
\address{Laboratoire J.-A. Dieudonn\'e\\
UMR 7351 du CNRS\\
Universit\'e de Nice\\
Parc Valrose\\
F-06108 Nice cedex 2, France}
\email{arnaud.beauville@unice.fr}
 
\begin{abstract}
We prove that any surface with $q=p_g=0$ embedded by a sufficiently large linear system 
admits a rank $2$ Ulrich bundle. In particular
every Enriques surface admits a rank $2$ Ulrich bundle.
\end{abstract}
\maketitle 
Let $X\subset \P^n$ be a complex projective variety of dimension $d$. A \emph{Ulrich bundle} on $X$ is a vector bundle $E$ on $X$ satisfying $H\pp(X,E(-1))=\ldots =H\pp(X,E(-d))=0$. This notion was introduced in \cite{ES}, where various other characterizations are given; let us just mention that it is equivalent to say that $E$ admits a linear resolution as a $\O_{\P^n}$-module, or that the pushforward of $E$ onto $\P^d$ by a general linear projection is a trivial bundle. 

In \cite{ES} the authors ask whether every projective variety admits a Ulrich bundle. The answer is known only in a few cases: hypersurfaces and complete intersections \cite{HUB}, del Pezzo surfaces \cite[Corollary 6.5]{ES}, abelian surfaces \cite{B}, sufficiently general K3 surfaces 
 \cite{AFO}. We show here that the   Lazarsfeld-Mukai construction produces a rank 2 Ulrich bundle on any surface with $q=p_g=0$, provided the linear system $|\O_S(1)|$ is large enough\,:
 \begin{thm*}
Let  $S\subset \P^{n}$ be a surface with $q=p_g=0$. Assume that the linear system $|K_S^{-1}(1)|$ contains an irreducible curve. Then $S$ admits a rank $2$ Ulrich bundle. 

In particular, every Enriques surface $S\subset \P^n$ admits a rank $2$ Ulrich bundle.
\end{thm*}
\pr  
Let $C$ be an irreducible curve in $|K_S^{-1}(1)|$; we have $\O_S(1)_{|C}=\omega _C$. Let $A$ be a general line bundle  of degree $g(C)+1$ on $C$. Then $h^0(A)=2$, $h^1(A)=0$, and $A$ is generated by its global sections. Let $E$ be the kernel of the evaluation map $H^0(C,A)\otimes \O_S\twoheadrightarrow A$. This is a rank 2 vector bundle on $S$. From the exact sequence 
\[0\rightarrow E\rightarrow H^0(C,A)\otimes \O_S\rightarrow A\rightarrow 0\]we get $H\pp(S,E)=0$. Since $\det E= \O_S(-C)$,   $E(1)$ is isomorphic to $ E^*\otimes K_S$, thus $H\pp(S,E(1))=0$. Therefore $E(2)$ is a Ulrich bundle for $S$.\qed

\bigskip	
\noindent\emph{Remarks}$.-$ 1)  It might be that every Enriques surface carries actually a Ulrich line bundle.  In \cite{BN} the authors conjecture that  Enriques surfaces with no $(-2)$-curves admit a Ulrich line bundle, and prove that conjecture in a significant case. 

\medskip	
2) Minimal surfaces with Kodaira dimension 0 fall into 4 classes, namely K3, Enriques, abelian and bielliptic surfaces. A bielliptic surface $S$ admits a finite \'etale covering $\pi :A\rightarrow S$ by an abelian surface $A$; if $L$ is a very ample line bundle on $S$ and $E$ a Ulrich bundle for $(A,\pi ^*L)$, then $\pi _*E$ is a Ulrich bundle for $(S,L)$. Thus all minimal surfaces with Kodaira dimension 0 admit a Ulrich bundle, except possibly some special K3 surfaces.

\medskip	
3) Let $S\subset \P^n$ be a surface with $q=p_g=0$, such that $\O_S(1)\cong K_S^r$, with $r\geq 3$. One can show that the linear system $|K_S^{r-1}|$ contains an irreducible curve, hence $S$ admits a rank $2$ Ulrich bundle.


\bigskip	
\begin{thebibliography}{X-X-X}
\bibitem[AFO]{AFO} M. Aprodu, G. Farkas, A. Ortega\,: \textsl{Minimal resolutions, Chow forms and Ulrich bundles on $K3$ surfaces}. J. Reine Angew. Math., to appear.
 Preprint  \texttt{arXiv:1212.6248}. 
\bibitem[B]{B} A. Beauville\,: \textsl{Ulrich bundles on abelian surfaces}. Proc. Amer. Math. Soc., to appear. Preprint \texttt{arXiv:1512.00992}.

\bibitem[BN]{BN} L. Borisov, H. Nuer\,: \textsl{Ulrich Bundles on Enriques Surfaces}.
Preprint  \texttt{arXiv:1606.01459}.

\bibitem[ES]{ES} D. Eisenbud, F.-O. Schreyer\,: \textsl{Resultants and Chow forms via exterior syzygies}. J. Amer. Math. Soc. \textbf{16}  (2003), no. 3, 537-579. 

\bibitem[HUB]{HUB}  J. Herzog, B. Ulrich, J.  Backelin\,: \textsl{Linear maximal Cohen-Macaulay modules over strict complete intersections}.
J. Pure Appl. Algebra \textbf{71}  (1991), no. 2-3, 187-202.


\end{thebibliography}
\end{document}